\documentclass[11pt,twoside]{article}

\usepackage{setspace}
\usepackage[margin=1.5in]{geometry}

\usepackage{amsthm,amsmath,verbatim,fancyhdr}
\usepackage{amssymb}
\usepackage[square,numbers]{natbib}
\usepackage[english]{babel}

\usepackage{newlfont}
\usepackage{mathrsfs}
\usepackage[mathscr]{euscript}
\usepackage{bbm}
\usepackage[arrow,matrix]{xy}
\usepackage[bookmarks,breaklinks]{hyperref}

\usepackage{amsthm}

\theoremstyle{plain}
\newtheorem{thm}{Theorem}[section]
\numberwithin{equation}{section}
\newtheorem{cor}[thm]{Corollary}
\newtheorem{lemma}[thm]{Lemma}
\newtheorem{prop}[thm]{Proposition}

\newtheorem{defn}[thm]{Definition}
\newtheorem{rem}[thm]{Remark}

\newcommand{\A}{\mathbb{A}}

\newcommand{\R}{\mathbb{R}}
\newcommand{\Z}{\mathbb{Z}}
\newcommand{\ZZ}{\mathrm{Z}}

\newcommand{\G}{\mathbb{G}}
\newcommand{\C}{\mathbb{C}}

\newcommand{\m}{\mathfrak{m}}

\newcommand{\uu}{\mathbf{u}}
\newcommand{\vv}{\mathbf{v}}

\newcommand{\eq}{\operatorname{eq}}  
  
\newcommand{\infl}{\operatorname{inf}}  
\newcommand{\Mor}{\operatorname{Mor}}
\newcommand{\Ker}{\operatorname{Ker}}
\renewcommand{\Im}{\operatorname{Im}}
\newcommand{\End}{\operatorname{End}}

\newcommand{\Fix}{\operatorname{Fix}}
\newcommand{\Aut}{\operatorname{Aut}}

\newcommand{\Gal}{\operatorname{Gal}}

\newcommand{\chara}{\operatorname{char}}
\newcommand{\wo}{\operatorname{\setminus}}
\newcommand{\dcl}{\operatorname{dcl}}
\newcommand{\acl}{\operatorname{acl}}
\newcommand{\Def}{\operatorname{Def}}

\newcommand{\GL}{\operatorname{GL}}
\newcommand{\PGL}{\operatorname{PGL}}
\newcommand{\Br}{\operatorname{Br}}
\newcommand{\et}{\mathrm{\acute{e}t}}
\newcommand{\set}[1]{\{\ #1\ \}}
\newcommand{\suchthat}[2]{\{\ #1\ \mid\ #2\ \}}

\newcommand{\bu}{\bullet}

\newcommand{\ul}{\underline}
\renewcommand{\bar}{\overline}

\def\Ind{\setbox0=\hbox{$x$}\kern\wd0\hbox to 0pt{\hss$\mid$\hss}\lower.9\ht0\hbox to 0pt{\hss$\smile$\hss}\kern\wd0}
\def\Notind{\setbox0=\hbox{$x$}\kern\wd0\hbox to 0pt{\mathchardef
\nn=12854\hss$\nn$\kern1.4\wd0\hss}\hbox to
0pt{\hss$\mid$\hss}\lower.9\ht0 \hbox to
0pt{\hss$\smile$\hss}\kern\wd0}

\begin{document}

\setlength{\headheight}{15pt}
\pagestyle{myheadings}
\markboth{{\normalsize Elimination of generalised imaginaries and
    Galois cohomology}}{{\small D. Sustretov}}

\begin{center}
{\LARGE\bf Elimination of generalised imaginaries and
    Galois cohomology}
\end{center}
\vspace{2ex}
\begin{center} 

  {\normalsize Dmitry Sustretov}\footnote{ The research leading to
    these results has received funding from the European Research
    Council under the European Union’s Seventh Framework Programme
    (FP7/2007- 2013)/ERC Grant Agreement No. 291111.  The author was
    also supported by Hill Foundation Scholarship when the work on the
    material present in the article has started.}
\end{center}
\vspace{4ex}

\bibliographystyle{abbrvnat}

\begin{minipage}[t]{0.85\linewidth}
\begin{abstract}
  The objective of this article is to characterise elimination of
  finite generalised imaginaries as defined in \citep{hru-groupoids}
  in terms of group cohomology. As an application, I consider series
  of Zariski geometries constructed \citep{hz,limit, quantum-zg} by
  Hrushovski and Zilber and indicate how their non-definability in
  algebraically closed fields is connected to eliminability of certain
  generalised imaginaries.
\end{abstract}

{
\small
\tableofcontents            % generate and include a table of contents
}
\end{minipage}

\section{Introduction}

Let $M$ be a strongly minimal non-locally modular structure. It has
been a long-standing conjecture of Zilber's \citep{zilber-icm} that
$M$ interprets an algebraically closed field (this statement is one of
the clauses of a statement widely known as Zilber's trichotomy
principle). It has been disproved by Hrushovski \citep{hrush-sm} but
then proved by Hrushovski and Zilber in the important case of Zariski
geometries \citep{hz}. The notion of a Zariski geometry is a natural
axiomatisation of the properties of Zariski topology on Cartesian
powers of algebraic varieties over an algebraically closed field, and
of compact complex manifolds (\citep{zilber-mtag,zg}). A Zariski
geometry is a structure $M$ endowed with a topology on every Cartesian
power of $M^n$ such that definable sets in $M$ are constructible sets
in the topology. Moreover, the structure must have a notion of
dimension (for example, Krull dimension), that satisfies certain
properties.

The main result of \citep{hz} asserts that a strongly minimal
non-locally modular Zariski geometry $M$ interprets an algebraically
closed field $k$. The next natural question to ask is that of the
possibility of co-ordinatisation of $M$, i.e. finding a definable in
$M$ injection $M \hookrightarrow k^n$ for some $n$, so that $M$
becomes a quasi-projective algebraic curve. Theorem~B of \citep{hz}
asserts that if one assumes just non-local modularity then one can
only guarantee the existence of a definable map $M \hookrightarrow
k^n$ with finite fibres, and in order to ensure that it is an
injection, one has to impose an extra condition: existence of a family
of one-dimensional sets in $M^2$ that separates points.

In Theorem~C of \citep{hz} the authors give an example of a Zariski
structure that projects onto an algebraic curve but the extra
structure on the fibres of the projection prevents this structure from
being interpretable in an algebraically closed field. 

Later, many more examples of Zariski geometries that are not
interpretable in an algebraically closed field have been constructed
by Zilber (\citep{limit, qho, quantum-zg}), not only in Morley rank
1. The work that lead to the results presented in this article started
as an attempt to find a uniform approach to proving
non-interpretability of such structures in an algebraically closed
field, identifying an obstruction that would appear in all known
examples. It was also desirable to be able to decide if the structures
in question were interpetable in a compact complex manifold considered
as a first-order structure.

This article suggests the following approach to this problem, building
on the notion of generalised imaginary sort introduced by Hrushovski
\citep{hru-groupoids}. In \citep{poizat} Poizat proposed a
model-theoretic generalisation of an absolute Galois group: the group
$\Gal(\acl(A)/\dcl(A))$ of automorphisms of $\acl(A)$ that fix
$\dcl(A)$ can be endowed with a topology generated by the base
consisting of the stabilisers of finite subsets of $\acl(A)$, and in a
theory that eliminates imaginaries there is a one-to-one
correspondence between closed subgroups of $\Gal(\acl(A)/\dcl(A))$ and
definably closed subsets of $\acl(A)$. That opens path for extending
the Galois cohomological results from theory of algebraically closed
fields to a general model-theoretic setting. Thus,
\mbox{Pillay}~\citep{pillay-galois} has noticed that the correspondence
between isomorphism classes of torsors and cocycle classes of the
first group cohomology group translates almost verbatim to
model-theoretic context.

In \citep{hru-groupoids} Hrushovski introduced generalised imaginaries
as certain sorts related to definable groupoids. Loosely speaking,
if one regards a definable groupoid as a generalisation of an
equivalence relation where equivalence classes can have automorphisms,
then generalised imaginary sort is something that is like an imaginary
sort, but also takes into account the automorphisms. 

In this article I introduce a notion of Morita equivalence of
definable groupoids (quite standard in other categories) so that
Morita equivalent groupoids give rise to bi-interpretable generalised
imaginary sorts. Groupoids that correspond to those imaginary sorts
which are interpretable in the home sort are called eliminable. The
notion of Morita equivalence gives the same equivalence classes of
definable groupoids as the notion of equivalence defined by Hrushovski
in \citep{hru-groupoids}. A generilised imaginary sort is a sort with
a structure of a definable groupoid torsor, i.e. a set acted upon
definably, freely and transitively by a definable groupoid.

I then prove (Theorem~\ref{cohomology}) that in a setting where
$K$-definable groupoids have a split torsor over $\acl(K)$, Morita
equivalence classes of connected $K$-definable groupoids that have an
Abelian isomorphism group of objects are in one-to-one correspondence
with classes in the second cohomology group of the absolute Galois
group of $K$. In this correspondence, eliminable groupoids correspond
to the trivial class.

% The proof relies on interpreting morphisms in the lower degree term
% long exact sequence of Hochschild-Serre spectral sequence (which is
% recalled in Appendix~\ref{hss}) in terms of operations on groupoids.

In Section~\ref{sec:brauer} I show that quantum Zariski geometries of
\citep{quantum-zg} are bi-interpretable with a certain generalised
imaginary sort. It is then easy to see that the definability of the
whole structure depends on eliminability of the corresponding
groupoid. Thus, eliminability can be decided by computing the
corresponding cohomology class. It also follows easily
(Section~\ref{sec:elim-ccm-omin}) that the interpretability of the
discussed structures in a compact complex manifold is equivalent to
interpretability in an algebraically closed field, since the
parameters of the groupoids corresponding to the obstruction sorts lie
in a projective variety, and thus their absolute Galois group is
exactly the same as in the algebraically closed fields.

The same situation is observed (Section~\ref{sec:groupext}) in the
case of non-standard Zariski structures defined in \citep{limit,hz},
if one assumes finiteness of the group action used in their
definition; there again the definability of the structure in an
algebraically closed field is equivalent to eliminability of a certain
generalised imaginary. 

I would like to thank Boris Zilber for introducing me to the circle of
problems that motivated the work presented in this article, and I am
grateful to Maxim Mornev, Martin Bays and Moshe Kamensky for many
helpful conversations.

\section{Generalised imaginaries}

\subsection{Groupoids and torsors}

Throughout the article we will write $X \times_{f,Z,g} Y$ for the fibre
product of definable sets, i.e. the set
$$
\suchthat{(x,y) \in X \times Y}{f(x)=g(y)}
$$
sometimes dropping $f$ and/or $g$ when they are clear from context.

\begin{defn}[Groupoid]
  A \emph{groupoid} is a category such that all its morphisms are
  isomorphisms. If a groupoid is small, i.e.\ if its objects and its
  morphisms are sets, then it is defined by the following data: a
  tuple $X_\bullet=(X_0, X_1)$ of sets along with maps $s,t,m,i,e$,
  where $s,t$ maps $X_1$ to $X_0$ (source and target objects), $c$
  maps $X_1 \times_{s,X_0,t} X_1$ to $X_1$ (composition of arrows),
  $i$ maps $X_1$ to itself (inverse), $e: X_0 \to X_1$, satisfying the
  natural axioms. 

  Let $\mathcal C$ be a category that has finite products. A
  \emph{groupoid $X_\bu$ internal to a category $\mathcal C$} is a
  pair of objects $X_0, X_1$ along with the morphisms $s,t,m,i,e$
  satisfying the mentioned identities.

  The set of morphisms from object $x$ to object $y$ is denoted
  $\Mor(x,y)$.  If $\Mor(x,x)$ is isomorphic to a group $A$ for all $x
  \in X_0$ then the groupoid $X_\bu$ is said to be \emph{bounded} by
  $A$.
\end{defn}

\begin{rem}
  The notation serves to underline the fact that a groupoid is in
  particular a simplicial set that only has $0$- and $1$-simplices.
\end{rem}

\begin{defn}[Definable groupoid]
  Let $\Def(\mathcal U)$ be the category of sets and maps definable
  with parameters in a monster model $\mathcal U$ of a complete theory
  $T$. Then a groupoid internal to $\Def(\mathcal U)$ is called
  \emph{definable groupoid} (cf. \citep{hru-groupoids}).

  A definable groupoid definable over a set of parameters $K$ is said
  to be \emph{bounded} by a definable group $A$ (definable over $K$ as
  well) if for any $x \in X_0$ there exists a definable isomorphism
  between $\Mor(x,x)$ and $A$, definable perhaps over a bigger set of
  parameters $L \subset K$.
\end{defn}

\begin{defn}[Action groupoid]
  Let $G$ be a group and let $\cdot: G \times X \to X$ be a group
  action. The \emph{action groupoid} is defined to be the groupoid
  with the morphisms $G \times X$ and objects $X$ where $s(g,x)=x$ and
  $t(g,x)=g\cdot x$, $(g,x) \cdot (h,gx)=(gh,x)$, and other structure
  maps defined in the obvious way.
\end{defn}

\begin{defn}[Groupoid quotient]
  Let $X_\bu$ be a groupoid. Let $E$ be the equivalence relation on
  $X_0$ which is the image of the map $(s,t): X_1 \to X_0 \times
  X_0$. The quotient $X_0/E$ is called the \emph{groupoid quotient}.
  We will denote it as $[X_\bu]$.
\end{defn}

\begin{defn}[Groupoid torsors]
\label{torsors}
  Let $X_\bu$ be a groupoid. A \emph{groupoid homogeneous space for
    $X_\bu$ over $Y$} is a map $p: P \to Y$ together with the
  \emph{anchor map} $a: P \to X_0$ and \emph{action map} $\cdot: X_1
  \times_{s,X_0,a} P \to P$ which commutes with the projection to
  $Y$. A homogeneous space is called \emph{principal} (or a
  \emph{torsor}) if for any two $f,g \in P$ such that $p(f)=p(g)$
  there exists a unique $m \in X_1$ such that $f \cdot m=g$.

  A morphism of groupoid torsors $P$ and $Q$ is a map $\alpha: P \to
  Q$ that respects the anchor map and commutes the action map:
  $a(\alpha(f))=a(f)$, $\alpha(m \cdot f)=m \cdot \alpha(f)$ for any
  $a \in X_\bu)$ and any $m \in \Mor(a,s(f))$.

  A groupoid $X_\bu$ is called \emph{eliminable} if there exists a
  $X_\bu$-groupoid torsor over $[X_\bu]$.
\end{defn}

Informally, a groupoid torsor is a collection of arrows from $Y$ to
$X_0$ with a possibility to compose these arrows with morphisms of
$X_\bu$ with the suitable source object. Note that if $X_\bu$ is a
groupoid with a single object $x$, a groupoid torsor is the same as
group $\Mor(x,x)$-torsor.

Let $E$ be an equivalence relation on a definable set $X_0$.
Transitivity, symmetry and reflexivity $E$ imply the existence of the
natural composition, inverse and identity structure maps of the
groupoid with the morphisms set $X_1=E \subset X_0 \times X_0$. Then
$X_0$ is a groupoid torsor over $X_0/E$. In a theory that does not
eliminate imaginaries $X_0/E$ lives in an imaginary sort. Elimination
of imaginaries is the condition that $X_0/E$ is in definable bijection
with a definable set in some of the home sorts of the theory.

In a similar vein, one might want to define generalised imaginary
sorts as sorts that contain a groupoid torsor for a groupoid (which
does not necessarily come from an equivalence relation). 

\begin{defn}[Generalised imaginary sort]
\label{gen-imag}
Consider a theory with elimination of imaginaries. Let $X_\bu$ be a
definable groupoid. A \emph{generalised imaginary sort} is an
expansion of the theory with an additional sort $P$ which is the
groupoid torsor for $X_\bu$ over $[X_\bu]$. The expansion includes the
maps $s: P \to X_0$, $t: P \to [X_\bu]$ and the action map $a: P
\times X_1 \to P$. The axioms that define $P$ to be a torsor are
clearly first-order.
\end{defn}

\subsection{Morita equivalence}

Suppose two groups $G, H$ act freely on two spaces $X, Y$ so that $X/G
\cong Y/H$. If one declares the corresponding action groupoids
equivalent then the set of action groupoids for free aciton modulo the
equivalence is the same as the set of groupoid quotients up to
isomorphism. One way to motivate the notion of Morita
equivalence  is that it generalises the equivalence just described to
arbitrary groupoids, in particular to arbitrary action groupoids,
taking into account the stabilisers. The following definitions and
lemmas quite standard, their analogues for the differential category,
for exapmple, can be found in \citep{behrend}.

\begin{defn}[Morita equivalent groupoids]
  A \emph{Morita morphism} $f_\bu: X_\bu \to Y_\bu$ is a pair of maps
  $f_0: X_0 \to Y_0, f_1: X_1 \to Y_1$ such that  the diagram
  $$
  \xymatrix{
    X_1 \ar[r] \ar[d]^{f_1} & X_0 \times X_0 \ar[d]^{f_0 \times f_0}\\
    Y_1 \ar[r] & Y_0 \times Y_0\\    
  }
  $$
  commutes, $f_0$ is surjective and for any $(x_1, x_2) \in X_0 \times
  X_0$ the map $f_1$ induces a bijection between $Mor(x,y)$ and
  $Mor(f_0(x),f_0(y))$. If one looks at groupoids as small categories,
  then the above conditions say precisely that Morita morphism
  defines a fully faithful functor which is surjective on objects.

  Two groupoids $X_\bu$ and $Y_\bu$ are called \emph{Morita
    equivalent} if there exists a third groupoid $Z_\bu$ together with
  two Morita morphisms $Z_\bu \to X_\bu$ and $Z_\bu \to Y_\bu$.
\end{defn}

This is indeed an equivalence relation: given equivalences\linebreak
\mbox{$X_\bu \leftarrow R_\bu \to Y_\bu$} and $Y_\bu \leftarrow S_\bu
\to Z_\bu$ consider the groupoid $R_1 \times_{Y_1} S_1
\rightrightarrows R_0 \times_{Y_0} S_0$, one checks that the naturally
defined Morita morphisms from this groupoid to $X_\bu$ and $Z_\bu$
define a Morita equivalence.

\begin{prop}
Let $X_\bu$ and $Y_\bu$ be Morita equivalent. Then $[X_\bu] \cong [Y_\bu]$.
\end{prop}

\begin{proof}
Indeed, one concludes easily from the definitions that if $f: Z_\bu \to
X_\bu$ is a morita morphism, then $a \sim_X b$ if and only if $f_0(a)
\sim_Y f(b)$ where $\sim_X, \sim_Y$ are the equivalence relations on
$X_0, Y_0$, images of $X_1, Y_1$ in $X_0 \times X_0, Y_1 \times Y_1$
respectively.
\end{proof}

\begin{cor}
  Connectedness is preserved under Morita equivalence.
\end{cor}

A class of Morita equivalent groupoids is morally a quotient space
``that remembers stabilisers''. 

\begin{lemma}
\label{tors-preimag}
  Let $f: X_\bu \to Y_\bu$ be a Morita morphism and let $P$ be a
  $Y_\bu$-torsor over $S$. Then there exists a $X_\bu$-torsor $f^{-1}(P)$ over $S$.
\end{lemma}

\begin{proof}
  Consider the torsor $X_0 \times_{f_0,Y_0,a} P \to S$ with the natural anchor map
  and the action map
  $$
  \sigma (x,p) = (x,f_1(\sigma)p)
  $$
  The commutative diagram in the definition of Morita morphism ensures
  that the action map commutes with a map to $X_0$.
\end{proof}

\begin{lemma}
\label{tors-imag}
Let $f: X_\bu \to Y_\bu$ be a Morita morphism and let $P$ be a
$X_\bu$-torsor over $S$. Then there exists a $Y_\bu$-torsor $f(P)$
over $S$.
\end{lemma}

\begin{proof}
  Consider the torsor $P /\sim\ \to S$ where $\sim$ is defined as
  $$
  \begin{array}{cll}
  p \sim q &\textrm{ if and only if }& \sigma p = q \textrm{ for some }
  \sigma \in \Mor(a(p),a(q)), \\
  & & \textrm{ such that } f_0(a(p))=f_0(a(q)) \textrm{ and
  }f_1(\sigma) = \mathrm{id}\\ 
  \end{array}
  $$
  The anchor map is the composition of the anchor map with $f_0$. The
  action map is defined as $\sigma [p] = [f^{-1}(\sigma)\cdot p]$.
\end{proof}

\begin{prop}
\label{morita-biinterp}
  Let $X_\bu$ and $Y_\bu$ be Morita equivalent definable
  groupoids. Then expansions of the theory with generalised imaginary
  sorts corresponding to $X_\bu$ and $Y_\bu$ are
  bi-interpretable. 
\end{prop}

\begin{proof}
  Let $X_\bu$ and $Y_\bu$ be Morita equivalent via Morita morphisms
  $f: Z_\bu \to X_\bu$ and $g: Z_\bu \to Y_\bu$. Let $P$ be the
  generalised imaginary sort associated to $X_\bu$. Then it follows
  from Lemmas~\ref{tors-preimag} and \ref{tors-imag} that the torsor
  $g(f^{-1}(P))$ is interpretable in the expansion of the theory with
  $P$ and is the generalised imaginary sort associated to $Y_\bu$.
\end{proof}

\begin{cor}
  Generalised imaginary sort that corresponds to an eliminable
  \hyphenation{group-oid}groupoid is interpretable in the base structure.
\end{cor}

\begin{prop}
\label{morita-bitors}
Let $X_\bu, Y_\bu$ be groupoids definable over $K$. Then $X_\bu$ is
Morita equivalent to $Y_\bu$ if and only if there exists a set $Q$
(\emph{$X_\bu$-$Y_\bu$-bitorsor}) which has the structure of a
 $X_\bu$-torsor over $Y_0$ and $Y_\bu$-torsor over $X_0$ such that the
 actions of $X_\bu$ and $Y_\bu$ commute.
\end{prop}

\begin{proof}    
Let $X_\bu \to Y_\bu$ be a Morita morphism. Then $Q = X_1
\times_{X_0} Y_0$ has a natural structure of a
$X_\bu$-$Y_\bu$-bitorsor. If $Q$ is a $X_\bu$-$Z_\bu$-bitorsor and $P$ is a
$Z_\bu$-$Y_\bu$-bitorsor then $(Q \times_{Z_0} P)/Z_1$ is a
$X_\bu$-$Y_\bu$-bitorsor.

Conversely, let $Q$ be a $X_\bu$-$Y_\bu$-bitorsor. Then $Y_1
\times_{Y_0} Q \times_{X_0} X_1 \rightrightarrows Q$ has a natural
structure of a groupoid, with natural Morita morphisms to $X_\bu$ and
$Y_\bu$. 
\end{proof}

\begin{cor}
  Let $X_\bu$ be a groupoid definable over $K$ and suppose $X_\bu$ is
  bounded by $A$. Then $X_\bu$ is eliminable if and only if it is Morita
  equivalent to the trivial action groupoid over $[X_\bu]$.
\end{cor}

Hrushovski defines equivalence of definable groupoids in
\citep{hru-groupoids} as follows: $X_\bu$ is equivalent to $Y_\bu$ if
there exist full definable functors $X_\bu \hookrightarrow Z_\bu,
Y_\bu \hookrightarrow Z_\bu$ that are injective on objects and such
that images of $X_\bu$, $Y_\bu$ in $Z_\bu$ meet every isomorphism
class in $Z_\bu$. Such a pair of functors gives rise naturally to a
$X_\bu$-$Y_\bu$-bitorsor of arrows from objects of $X_\bu$ to objects
of $Y_\bu$ therefore, by Proposition~\ref{morita-bitors}, $X_\bu$ and
$Y_\bu$ are Morita equivalent. The inverse implication is slightly
less straightforward.

\begin{prop}
\label{hrush-morita}
Let $X_\bu \leftarrow Z_\bu \to Y_\bu$ be a Morita equivalence. Then
there exist definable full faithful functors from $X_\bu$ and $Y_\bu$
into a groupoid $W_\bu$ which are injective on objects and such that
the images meet every isomorphism class of objects of $W_\bu$.
\end{prop}

\begin{proof}
  Let $Q$ be the $X_\bu$-$Y_\bu$-bitorsor witnessing the Morita
  equivalence of $X_\bu$ and $Y_\bu$ (Lemma~\ref{morita-bitors}). 

  Define $W_0 = X_0 \sqcup Y_0, W_1 = X_1 \sqcup Q \sqcup Q^{-1}
  \sqcup Y_1$ where $Q^{-1}$ is a copy of $Q$. There are natural
  source and target maps defined on $Q, Q^{-1}$ coming from bitorsor
  structure. Thinking of elements of $Q$ as arrows from $X$ to $Y$
  and of elements of $Q^{-1}$ as arrows from $Y$ to $X$ define
  composition of arrows: the composition of arrows in $X_1$ and $Y_1$
  are given by composition in respective groupoids, the composition of
  arrows in $Y_1$ and $X_1$ and arrows in $Q, Q^{-1}$ is given by
  bitorsor structure.

  The fact that inclusions $X_\bu \hookrightarrow W_\bu, Y_\bu
  \hookrightarrow W_\bu$ are fully faithful follows from the fact that
  action of $X_1$ and $Y_1$ on $Q, Q^{-1}$ is free. 
\end{proof}

The following statement generalises the lemma of Lascar and Pillay
about elimination of imaginaries up to finite ones in strongly minimal
theories.

\begin{prop}
\label{finite-quot}
  Let $M$ be strongly minimal and let $\acl(\emptyset) \cap M \neq
  \emptyset$. Let $X_\bu$ be a groupoid with $X_0 \subset M^n$. Then
  there exists a Morita equivalent groupoid $Y_\bu$ such that $Y_0 \to
  [Y_\bu]$ is finite. 
\end{prop}

\begin{proof}
  Let $R$ be the image of $X_1$ in $X_0 \times X_0$. By
  Lemma~1.6~\citep{pillay-acf}, there exists a definable $Y_0 \subset
  X_0$ such that the $R$-equivalence classes of $Y_0$ are finite. The
  groupoid which is the restriction of $X_\bu$ to $Y_\bu$ clearly
  satisfies the requirements of Propostion~\ref{hrush-morita} and
  hence is Morita equivalent to $X_\bu$.
\end{proof}

Finally, let us remark that the the notion of \emph{retractability}
defined in \cite{goodrick} is equivalent to having a Morita morphism
to a groupoid such that $X_0=[X_\bu]$; one uses the straightforward
generalisation of the argument in Proposition~1.11, \emph{loc.cit.},
applying it to groupoids that are not necessarily connected.

\section{Galois cohomology}

I will work in a theory that eliminates imaginaries, so the
Galois correspondence for strong types applies:

\begin{thm}[\citet{poizat}, Theorem~14]
  Let $K \subset L \subset \acl(K)$, then $\Aut_L(\acl(K))$ is the
  closed subgroup of $\Aut_L(\acl(K))$ that fixes $\dcl(L)$.
\end{thm}

For any set $K$ denote $G_K:=\Aut(\acl(K)/\dcl(K))$. If $G_L \subset
G_K$ is normal I will denote the quotient $\Gal(L/K)$. If $A$ is a
definable group then I denote $A(K)$ the group of tuples definable
over $K$ that belong to $A$; I will also denote the algebraic closure
of a set $K$ as $\bar{K}$. If $f: X \to Y$ is a map between sets
defined over $K$ and $\sigma$ is an element of $G_K$, then $\sigma(f):
\sigma(X) \to \sigma(Y)$ will be the map obtained from $f$ by
conjugation by $\sigma$: $\sigma(f) = \sigma^{-1} f \sigma$.

For a connected groupoid $X_\bu$ let us say that a $X_\bu$-torsor $P$
is \emph{split over $L$} if it has a point $x \in P$ definable over
$L$.

\begin{lemma}
  \label{torsor-aut}
  Let $P$ be a torsor of a connected groupoid $X_\bu$ bounded by a
  definable group $A$ and split over $L$. The group of definable (over
  some set $M \supset L$) automorphisms of $P$ is in bijective
  correspondence with elements of the group $A(M)$.
\end{lemma}

\begin{proof}
  Straightforward.
\end{proof}

\begin{thm}
\label{cohomology}
Let $K$ be a set of parameters. Let $A$ be a definable over $K$
Abelian group.  There exists a bijective correspondence between Morita
equivalence classes of connected groupoids definable over $K$,
eliminable over $\acl(K)$, with torsors split over $\acl(K)$ and
bounded by $A$ and cohomology classes in $H^2(G_K, A(\bar K))$.
% $$
% \left\{
% \begin{array}{c}
% \textrm{Morita equivalence classes of }\\
% \textrm{connected groupoids }\\
% \textrm{definable  over } K,\\
% \textrm{eliminable over } \acl(K),\\
% \textrm{with torsors split over } \acl(K),\\
%  \textrm{ and bounded by } A\\
% \end{array}
% \right\}
% \Leftrightarrow
% \left\{
% \begin{array}{c}
% \textrm{cohomology classes in } H^2(G_K, A(\bar K))
% \end{array}
% \right\}
% $$
Eliminable groupoids correspond to the trivial cohomology class. 
\end{thm}

\begin{proof}
  $ $\\
  \ul{Step 1}: describe a map that sends a cocycle that reperesents a
  class in $H^2(G_K,A(K))$
  to a groupoid with  a torsor over $\bar K$.\\
  Let $X_\bu$ be a groupoid and let $P$ be an $X_\bu$-torsor definable
  over a set of parameters $L \supset K$. If $L$ is minimal such that
  $G_L$ is normal in $G_K$ then the orbit of $P$ under the action of
  $G_K$ consists of several copies of $P$ that are in one-to-one
  correspondence with elements of $\Gal(L/K)$. I will denote the
  Galois conjugates of $P$ as $P_\sigma$ for all $\sigma \in G_K$ even
  though by doing so I denote a particular Galois conjugate by several
  different names. Until the end of this proof the composition will be
  written without the symbol $\circ$ and from left to right.

  Pick a point $x \in P$ definable over $\bar K$. Choose a continuous
  section $j: \Gal(L/K) \to G_K$. Let $v_\sigma: P \to P_\sigma$
  be the isomorphism of $X_\bu$-torsors that sends $x$ to
  $\sigma(x)$. Let $Q=\bigcup\limits_{\sigma \in G_K} P_\sigma$
  and define $u_\sigma: Q \to Q$ as follows:
  $$
  u_\sigma (y)  = v_{j(\alpha)^{-1}} u_\sigma v_{j(\alpha)}
  $$
  for $y \in P_\alpha$. The maps $u_{\sigma}$ can be thought of as the
  \emph{scindage} in the terminology of \citep{giraud}, IV.3.5. Note
  that $u_\sigma f = \sigma (f) u_\sigma$ for any $f: Q \to Q$ any $\sigma \in
  G_K$.

  Define
  $$
  h(\sigma,\tau) = u_\sigma u_\tau u_{\sigma\tau}^{-1}
  $$  
  The expression on the right is an automorphism of $P$ and therefore
  can be identified with an element of $A(\bar K)$ by
  Lemma~\ref{torsor-aut}.

  Let us check that this indeed defines a cocycle, i.e. that the
  equality
  $$
  \alpha(h(\sigma,\tau)) h(\alpha, \sigma\tau) =
  h(\alpha\sigma, \tau) h(\alpha, \sigma)
  $$
  holds. Indeed,
  $$
  u_\alpha u_\sigma u_\tau =  h(\alpha, \sigma)
  u_{\alpha\sigma} u_\tau = h(\alpha,\sigma)
  h(\alpha\sigma,\tau) u_{\alpha\sigma\tau}
  $$
  but on the other hand
  $$
  u_\alpha u_\sigma u_\tau = u_\alpha h(\sigma,\tau)
  u_{\sigma\tau}=h(\sigma,tau) u_\alpha u_{\sigma\tau} =
  \alpha(h(\sigma,\tau)) h(\alpha,\sigma\tau) u_{\alpha\sigma\tau}
  $$

  \ul{Step 2}: describe a map that puts a groupoid and a torsor into
  correspondence to a cocycle in $H^2(G_K, A(\bar K))$.\\
  By Lemma~\ref{coho-limit} for any $\alpha \in H^2(G_K, A(\bar K))$
  there exists $\beta \in H^2(\Gal(L/K),\linebreak A(L))$ for some $L
  \subset \bar K$ such that $\alpha=\infl(\beta)$. Let $G$ be the
  definable group, an extension of $\Gal(L/K)$ by $A$, that
  corresponds to (Theorem~\ref{group-ext}) the cocycle
  $\beta$. Consider the action of $G$ on a $\Gal(L/K)$-torsor $P$
  where an element $g \in G$ acts as its projection $p(g) \in
  \Gal(L/K)$, and let $X_\bu$ be the associated action groupoid. Then
  $\cup \Mor (y,x)$ for $x \in X_0$ is definable over $L$ (it is a
  union of $|\Gal(L/K)|$ copies of $A$) and is naturally an
  $X_\bu$-torsor over a singleton.

  \ul{Step 3}: show that groupoids corresponding to cohomologous
  cocycles are  Morita-equivalent.\\
  This amounts to showing (by Step~2 and Theorem~\ref{group-ext})
  that if $f: G \to G'$ is an isomorphism of group extensions, i.e. if
  $$
  \xymatrix{
  1  \ar[r] & A  \ar@{=}[d] \ar[r]   & G \ar[d]^f  \ar[r] & \Gal(L/K)
  \ar@{=}[d] \ar[r] & 1 \\  
  1  \ar[r] & A   \ar[r] & G'  \ar[r] & \Gal(L/K) \ar[r]  & 1 \\
  }
  $$
  is a commutative diagram then the action groupoids for the actions
  of $G$ and $G'$ on a $\Gal(L/K)$-torsor are Morita equivalent. In
  fact, it is clear by construction of the action groupoid that they
  are definably isomorphic.

  \ul{Step 4}: show that cocycles corresponding to a groupoid $X_\bu$
  and two different $X_\bu$-torsors $P,P'$ are cohomologous.

  Pick some tuples $x \in P, y \in P'$ defined over $\bar K$ such that
  $a(x) = a(y)$ (where $a$ is the anchor map).  We will identify $P$
  and $P'$ via an isomorphism that sends $x$ to $y$. Let $\eta: P' \to
  P$ the isomorphism of torsors that sends $x$ to $y$.  Let $u_\sigma,
  u'_\sigma$ be the maps as in Step~1 used to obtain the cocycles
  $h,h'$ for $P$ and $P'$ respectively. Define a cochain $g: G_K \to
  A(\bar K)$:
  $$
  g(\sigma) = u'_\sigma u_\sigma^{-1}
  $$
  (here $g(\sigma)$ is identified with an element of $A(\bar K)$ by
  Lemma~\ref{torsor-aut}). Let $h, h' \in H^2(G_K, A(\bar K)$ be
  cocycles that are obtained using procedure from Step~1 from torsors
  $P$ and $Q$. Then
  $$
  \begin{array}{ll}
  h'(\sigma, \tau) & = g(\sigma) u_\sigma g(\tau) u_\tau
  u_{\sigma\tau}^{-1} g(\sigma\tau)^{-1} = \\
  & = g(\sigma) u_\sigma g(\tau) u_\sigma^{-1} u_\sigma u_\tau
  u_{\sigma\tau}^{-1} g(\sigma\tau)^{-1} = \\
  & = g(\sigma)  \sigma(g(\tau))   h(\sigma, \tau)  g(\sigma\tau)^{-1}
  = \\
  & = g(\sigma)  \sigma(g(\tau)) g(\sigma\tau)^{-1}   h(\sigma, \tau)  \\
   \end{array}
  $$  
  and therefore $h$ and $h'$ are cohomologous.

  Consequently, if a groupoid has a torsor definable over $K$ then the
  associated cocycle is cohomologous to zero. Indeed, in the latter
  case, as all the maps $u_\sigma$ of Step~1 are automorphisms of the
  $K$-definable torsor $P$, and $u_\sigma=\sigma$, therefore the
  corresponding cocycle is the zero cocycle.
  
  \ul{Step 5}: show that cocycles corresponding to Morita equivalent
  groupoids are the same.

  Let $Y_\bu$ have a torsor $Q$ and let $f: X_\bu \to Y_\bu$ be a Morita
  morphism. Then by Lemma~\ref{tors-preimag} $P= X_0
  \times_{f_0,Y_0,a} Q$ is a $X_\bu$-torsor. Let $u_\sigma: Q \to
  Q_\sigma$ be a collection of isomorphisms of torsors, then
  $$
  v_\sigma (x,q) = (x,u_\sigma(q))
  $$
  is a collection of isomorphisms of Galois conjugates of $Q$.  It is
  easily checked that $v_\sigma$ gives rise to exactly the same
  cocycle.

  Similarly, if $X_\bu$ has a torsor $P$ and $f: X_\bu \to Y_\bu$ is a
  morphism of torsors then $Q=P/\sim$ is a $Y_\bu$-torsor where $\sim$
  is the equivalence relation defined in Lemma~\ref{tors-imag}. The
  morphisms $u_\sigma$ are determined by the choice of a point $x \in
  P$ and they descend to morphisms $v_\sigma$ between Galois
  conjugates of $Q$ determined by the equivalence class $x/\sim \in
  Q$. One checks that again $v_\sigma$ give rise to the same cocycle
  as $u_\sigma$.  
  
  \ul{Step 6}: show that the correspondence is a bijection.

  Let $P$ be a torsor of a groupoid $X_\bu$. If $P$ is definable over
  a set of parameters $L$ such that $G_L$ is normal in $G_K$ then the
  cocycle $h$ representing an element in $H^2(G_K,A(\bar K))$
  constructed in Step~1 is the inflation of a cocycle $\eta$
  representing an element in $H^2(\Gal(L/K), A(L))$. Let $G$ be the
  definable group extension of $\Gal(L/K)$ by $A$ corresponding to
  $\eta$. Let us show that the group action groupoid for $G$ acting on
  a $\Gal(L/K)$-torsor via projection $G \to \Gal(L/K)$, call it
  $Y_\bu$, is Morita equivalent to $X_\bu$.

  Let $Q=\bigsqcup\limits_{\sigma \in \Gal(L/K)} P_\sigma$ and let $B$
  be the $\Gal(L/K)$-orbit of the tuple of elements (defined over $L$)
  used to define $P$. Then, by Galois invariance, there is a
  projection $p:Q \to B$ defined over $K$.  Let $Z_\bu$ be the
  groupoids such that $Z_0=X_0 \sqcup Y_0$ and $Z_1=X_1 \sqcup Y_1
  \sqcup Q \times Q'$ where $Q'$ is a copy of $Q$. The composition of
  arrows in $X_1$ and in $Y_1$ is defined as in groupoids $X_\bu$ and
  $Y_\bu$, the elements of $X_1$ compose with elements of $Q$ as
  follows:
  $$
  f \cdot \sigma(x) = \sigma (f \cdot x)
  $$
  the elements of $Q$ compose with elements of $Y_1$ as follows
  $$
  x \cdot g = u_{p(g)}(x)
  $$
  and similarly for inverse arrows. 

  Now by Lemma~\ref{hrush-morita} $X_\bu$ and $Y_\bu$ are Morita
  equivalent.

  Let $h$ be a cocycle representing an element of $H^2(G_K,A(\bar K))$
  that is an image of an element of $\eta \in H^2(\Gal(L/K),A(L))$,
  let $X_\bu$ be the action groupoid corresponding to the extension
  described by the element $\eta$. It is straightforward to check, by
  following the construction of the Step~1, that the cocycle that
  corresponds to $X_\bu$ is $\eta$.
\end{proof}

\begin{rem}
  If $\acl(K)$ is a model then the theorem allows to classify all
  connected groupoids bounded by an Abelian group definable over $K$. 
\end{rem}

\begin{rem}
  \label{remark-any-eq-action} 
  Note that as a part of the proof (Step~6) we have seen that a
  connected groupoid that has a torsor, defined over $\acl(K)$, is
  Morita equivalent to an action groupoid for an extension $G$ of the
  group $\Gal(L/K)$ for some $L \supset K$ by an Abelian definable
  group $A$, acting on a $\Gal(L/K)$-torsor via the projection on
  $\Gal(L/K)$.
\end{rem}

\begin{prop}
  \label{group-image}
  Let $f: A \to B$ be a definable map of definable groups and let
  $\alpha \in H^2(G_K, A(\bar K))$. Then the generalised imaginary
  sort corresponding to $f^*(\alpha) \in H^2(G_K, B(\bar K))$ is
  interpretale in the generalised imaginary sort corresponding to
  $\alpha$.
\end{prop}

\begin{proof}
  Let $X_\bu$ be the groupoid definable over $K$ corresponding to the
  cocycle $\alpha$ with a torsor $P$ and let $Y_\bu$ be the groupoid
  corresponding to the cocycle $f^*(\alpha)$. We may assume that
  $X_\bu$ and $Y_\bu$ are group action torsors, so in particular the
  action of $A$ on $P$ is defined. Then $Y_\bu$ has the torsor $Q=(P
  \times B)/A$ where the action of $A$ is defined as
  $$
  a \cdot (p,b) = (a^{-1} \cdot p, f(a)b)
  $$
  One checks that $Q$ is a $Y_\bu$-torsor.
\end{proof}

% \begin{rem}
%   The discussion at the end of Section~3 in \citep{hru-groupoids} can
%   be reinterpreted in terms of well-known results about group
%   cohomology. Indeed, $H^2(\Z/2\Z, A)$ is known to be non-trivial
%   unless $A$ is 2-torsion, therefore one finds plenty of non-split
%   covers of $\ACF_L$, where the absolute Galois group of $L$ is
%   $\Z/2\Z$. The eliminability of finite Abelian groupoids over PAC
%   fields is the consequence of the fact that the Brauer group of a PAC
%   field is trivial (\citep{fried-jarden}).
  
%   The correspondence given in Theorem~\ref{gerbescohom} seems to
%   provide an answer to the question raised in \citep{hru-groupoids} (
%   Problem~3.11). Let $K$ be a finitely generated field. Any extension
%   of the absolute Galois group of $K$ gives rise to a class in
%   $H^2(G_K,A)$ and hence to a gerbe and a generalised imaginary sort.
% \end{rem}

\section{Non-standard Zariski structures}

In this section I will look at two series of examples of Zariski
structures constructed by Hrushovski and Zilber, and examine their
interpretability in various theries. It will turn out that
intrepretability is closely related to eliminability of certain
generalised imaginaries.

\subsection{Group extensions}
\label{sec:groupext}

Let $X$ be an algebraic variety defined over an algebraically closed
field $k$. Let $H$ be an abstract group that acts on $X(k)$ such that
the stabilizer of any point is either $G$ or the trivial
subgroup. 

Consider some extension $1 \to A \to G \to H \xrightarrow{\pi} 1 $ with $A$
finite. Let $\Fix(H)$ denote the set of fixed points of $H$. Each
orbit of $H$ on $X \wo \Fix(H)$ is a principal homogeneous space. Pick
a representative $x_\alpha$ in each orbit and let $D$ be the set
$$
\bigsqcup (G \cdot x_{\alpha} ) \sqcup \Fix(H)
$$
Let $p: D \to X$ be a map that maps elements of the form $g\cdot
x_\alpha$ to $\pi(g) \cdot x_\alpha$ and is identity on
$\Fix(H)$. The group $G$ acts naturally on $D$. If one declares the
pre-images of closed subsets of $D^n$ under $p$ and graphs of the
action by elements of $G$ closed, this defines a topology on $D$ which
satisfies the axioms of a Zariski geometry, as proved in
\citep[Proposition~10.1]{hz}. The Proposition also asserts that the
isomorphism type of the structure $D$ does not depend on the choices
of representatives in the orbits of $H$. Let us denote this structure
$D(X,G)$.  Hrushovski and Zilber further prove (Theorem~C) that there
exists a variety $X$ and a group extension such that $D(X,G)$ is not
interpretable in an algebraically closed field. More examples of
structures of the form $D(X,G)$ were considered in \citep{limit}.

If one assumes that $H,G$ are finite (and hence definable) it is easy
to see that construction is equivalent to adding a generalised
imaginary sort for an action groupoid. Indeed, define $X_\bu$ to be
the action groupoid for the lifting of the action of $H$ on $X$ to
$G$. Then $D$ is a $X_\bu$-torsor, with anchor map $p$:
$$
\xymatrix{
  X_1 = & X \times G \ar@<-0.5ex>[d]_s  \ar@<0.5ex>[d]^t
  & D \ar[dl]^p \ar[d]^{q \circ p} \\
  X_0 = & X & X/H\\
}
$$
where $q: X \to X/H$ is the natural projection map. The action of
$X_\bu$ is given by the action of $G$ on $D$.

\begin{prop}
  Let $G$ be finite. The structure $D(X,G)$ is interpretable in an
  algebraically closed field if and only if $X_\bu$ is eliminable.
\end{prop}
\begin{proof}
  The left to right direction tautologically follows from the
  definition of the structure.

  Let us prove the right to left direction. As follows from trichotomy
  for Zariski geometries, $D(X,G)$ interprets an algebraically closed
  field $k \subset D^n$, and moreover, there exists a definable
  embedding $X \to k^m$ (for 1-dimensional $X$, see \citep[Sections~6
  and 9]{hz}, for $X$ arbitrary, see \cite[Section~4.4]{zg}). By a
  theorem of Poizat (\citep[Theorem~4.15]{poizat-groups}), if $D(X,G)$
  is interpretable in an algebraically closed field $F$ then $k$ is
  definably isomorphic to $F$. Therefore, one may assume that $D$ is a
  definable set in $k$ and the projection $p$ is definable. But that
  means precisely that $X_\bu$ is eliminable.
\end{proof}

\begin{cor}
  If $X$ is a curve then any structure $D(X,G)$ with $G$ finite is
  interpretable in an algebraically closed field.
\end{cor}

\begin{proof}
  As follows from Tsen's theorem, $k(X)$ is quasi-algebraic and Galois
  cohomology of quasi-algebraic fields with coeffecients in torsion
  modules vanishes beyond degree 1 (\citep{serre-galois},II.\S 3.2),
  so $H^2(k(X), \mu_n)=0$. It follows that a
  restriction of $X_\bu$ to some subset $Y \subset X$ is
  eliminable. The restriction of $X_\bu$ to $X \wo Y$ is a groupoid
  with finitely many objects and hence is eliminable, perhaps after
  adding some parameters from $\acl(\emptyset)$.
\end{proof}

\subsection{Quantum Zariski geometries}
\label{sec:quantum}

The paper \citep{quantum-zg} considers a large class of Zariski
geometries that are constructed from a certain piece of data that (as
will be shown below) defines an Azumaya algebra over a variety. We now
recall the definition of a quantum Zariski geometry to fix notation
(Section~2, \citep{quantum-zg}).

\begin{defn}[Data for quantum Zariski geometry]
\label{quantum-zg}
  Let $k$ be an algebraically closed field, and $A$ be an associative
  unital finitely. The input data for a quantum Zariski geometry is
  \begin{enumerate}
  \item an associative algebra $A$ finitely generated over its center $\ZZ(A)$,
    and the center is a commutative algebra of finite type over $k$
    which is the coordinate ring of a variety $X$;
  \item a collection of irreducible modules $M_x$, all of fixed
    dimension $n$ over $k$, where $x$ ranges in the maximal spectrum
    of $\ZZ(A)$, such that $M_x$ is annihilated by $x$;
  \item\label{canbase} a choice of generators of $A$, $\uu_1, \ldots,
    \uu_d$, a choice of bases $e^{\alpha}_i(x) \in M_x, 1 \leq d \leq
    n$ (called \emph{canonical bases}), where $\alpha$ ranges in some
    finite set $B$, and a system of polynomials $\set{f_l(t,x)}, t
    = \{t_{ij}^k\}, 1 \leq i,j \leq n, 1 \leq k \leq d, x \in X$ such
     that for any $x \in X$ and any $t_{ij}^k$ that satisfy the
     equations $f_l(t,x)=0$  there is $\alpha \in B$ such that
     $\uu_i$ have the form
    $$
    \uu_k\, e^{\alpha}_j(x) = \sum_{i=1}^n\ t_{ij}^k\, e^{\alpha}_i(x)
    $$
  \item a finite group $\Gamma$ and a partial map $g: X \times \Gamma
    \to \GL_n(k)$ such that for all $\gamma \in \Gamma$, the map $g(-,
    \gamma) \to \GL_n(k)$ is regular on some open subset of $X$ for
    every $\gamma$,  and for any $x \in X$, $g(x, -)$ is defined on a
    subgroup of $\Gamma$ and is injective. For any $\alpha, \beta \in
    B$, there exists $\lambda \in k^\times$ and $\gamma \in \Gamma$
    $$
    e^{\beta}(x) = \lambda \sum_{j=1}^n g_{ij} (x,\gamma) e^{\alpha}(x)
    $$
  \end{enumerate}
  The equations, the map $g$ and the defining relations of the algebra
  $A$ can be defined over a subfield $k_0 \subset k$.
\end{defn}

\begin{rem}
  Note that the clause~4 effectively says that $g$ is a 1-cocycle of
  group cohomology of $\Gamma$ with values in the projective group
  $\PGL_n(k[U])$ for some open $\Gamma$-invariant $U \subset X$ with
  trivial action of $\Gamma$. In fact, nothing prevents one from
  considering a cocycle with values in $\PGL_n(k[Y])$ where $Y$ is the
  Galois cover of $U$ with Galois group $\Gamma$ defined by equations
  $f(x,t)=0$ from clause~3, so this is further allowed. 
\end{rem}

\begin{defn}[Quantum Zariski geometry]
\label{qzg-defn}
A \emph{quantum Zariski geometry} associated to the data described in
the previous definition is a structure with two sorts $V,k$ and a
projection map $p: V \to X(k)$ where
  \begin{enumerate}
  \item $k$ has the structure of an algebraically closed field, and
    $X$ is an affine variety with coordinate ring $Z(A)$, viewed as
    a definable subset in $k^m$;
  \item $V$ has a fibrewise structure of a $k$-vector space, i.e.\ the
    language on $V$ has graphs of operations $V \times V \to V$ and $k
    \times V \to V$ that restrict to graphs of addition and
    multiplication by a scalar on every fibre $V_x=p^{-1}(x), x \in X$;
  \item the structure contains graphs of maps $\uu_i: V \to V$ that
    restrict to linear self-maps on fibres of $p$, and are defined as
    follows. For every $x_0 \in X$, for every solution $t^k_{ij}$ of
    the system of equations $f_l(t,x_0)=0$ there exists a basis in
    $V_{x_0}$ such that $\uu_k$ acts on $V_{x_0}$ by the matrix
    $(t^k_{ij})$ in this basis.
  \end{enumerate}  
  Every fibre $V_x$ is therefore isomorphic to $M_x$ as an
  $A$-module.
\end{defn}

Zilber has shown in \citep{quantum-zg}, Lemma~2.4, that the structure
is unique up to isomorphism, once the base algebraically closed field
is fixed.

I am going to show that the described structure is bi-interpretable
with an algebraically closed field with an added generalised imaginary
sort. In order to do that I will argue that a quantum Zariski
geometry encodes an Azumaya algebra over the variety $X$.

\begin{defn}[Azumaya algebra]
  Let $A$ be an algebra over the ring of regular functions $k[U]$ of
  an affine variety $U$. The algebra $A$ is called \emph{twisted form
    of a matrix algebra} if for some morphism of varieties $Y \to X$
  the base change $A \otimes_{k[U]} k[Y]$ is isomorphic to the matrix
  algebra $M_n(k[Y])$. An algebra $A$ over $k[X]$ is called an
  \emph{Azumaya algebra} if for any point $x\in X$ there exist an open
  affine neighbourhood $U$ such that $A \otimes_{k[X]} k[U]$ is a
  twisted form of a matrix algebra. The rank of an Azumaya algebra is
  its rank as a module.
\end{defn}

From now on I adopt a simplifying assumption that the morphism $Y \to
X$ is Galois and the function $g: \Gamma \times X \to GL_n$ is regular
everywhere on $\Gamma \times X$. I will show that then the input data
of a quantum Zariski geometry defines two objects: a twisted matrix
algebra over $X$ and an injective morphism of $A$ into this algebra.

Indeed, the input data specifies a Galois cover $Y \to X$, where $Y$
is the subvariety of $X \times \A^{n^2\cdot d}$ defined by the
equations $f_l(x,t) = 0$ (clause~3 of the Definition~\ref{quantum-zg}),
and a 1-cocycle $g: \Gamma \to \PGL_n(k[Y])$ with values in the
group of regular functions from $Y$ to $\PGL_n$
(clause~4). Consider the matrix algebra $M_n(k[Y])$. The group
$\Gamma$ acts on $M_n(k[Y])$ by conjugation by $g(\gamma)$ for $\gamma
\in \Gamma$. Let $B$ be the $k[X]$-algebra $M_n(k[Y])^\Gamma$ of
elements fixed by this action.

By Galois descent (\citep{knusoj}, II.\S 5) there exists an
isomorphism $\iota: B \otimes_{k[X]} k[Y] \to M_n(k[Y])$ such that for
any $\gamma \in \Gamma$, $\iota^{-1} \circ \gamma \circ \iota =
g(\gamma)$.  As the equations $f_l(x,t)=0$ define $Y$, the variables
$t_{ij}^k$ are regular functions on $Y$, and the matrices $(t_{ij}^k)
\in M_n(k[Y])$ are part of the input data. Since these matrices are
$\Gamma$-invariant, they descend to elements of $B$. This defines a
map $\eta: A \to B$. 

Note that for a closed point $x\in X$, $B \otimes k(x) \cong
M_n(k(x))$ since $B \otimes k(x)$ is a central simple algebra over the
residue field at a point $k(x)$ (Proposition~IV.2.1, \citep{milne})
and $k(x)$ is algebraically closed, though the isomorphism is not
canonical.

\begin{prop}
\label{qzg-inj}
  The map $\eta$ is injective.
\end{prop}

\begin{proof}
  It follows in a straightforward way from the definitions that the
  $A$-module given by compostion of $\eta$ and the reduction map $B
  \to B\otimes k(x) \cong M_n(k)$ is isomorphic to $M_x$.

  Since the annihilator of $M_x$ is the maximal ideal $\m_x \subset
  k[X]$, the kernel of $\eta$ therefore is contained in the
  intersection of all maximal ideals of $k[X]$ which is zero since
  $k[X]$ is reduced.
\end{proof}

A classical fact from linear algebra (for proof see, for example,
\citep{burnside}) helps establish that $\eta$ is surjective if the
fibre modules are irreducible.

\begin{prop}[Burnside's theorem on matrix algebras]
 Let $V$ be a vector space over a field $k$, and let $A
 \subset \End(V)$. Suppose that $A$ does not have invariant
 subspaces. Then $A \cong M_n(V)$.
\end{prop}

\begin{prop}
\label{qzg-surj}
If all modules $M_x$ are irreducible then the image of the morphism
$\iota$ from the statement of Proposition~\ref{qzg-inj} coincides with $B$.
\end{prop}

\begin{proof}
  By Burnside's theorem the maps $\eta \otimes k(x): A \otimes k(x)
  \to B \otimes k(x)$ are isomorphisms for all $x \in X$. Therefore
  $A$ is an Azumaya algebra of the same rank as $B$ (see, for example,
  Proposition~III.2.1 in \citep{milne}). Since $f$ is an inclusion, it
  is an isomorphism.
\end{proof}

\subsection{Definability and Brauer group}
\label{sec:brauer}

There is a natural groupoid that one can assosciate to a quantum
Zariski geometry.

\begin{defn}[Splitting groupoid]
Let us keep notation for the data that defines a quantum Zariski
geometry from Definition~\ref{qzg-defn}. 

The splitting groupoid $S_\bu$ has the objects set $S_0=Y$ and the
morphisms set $S_1$ is $(Y \times_X Y) \times k^\times$ with the
obvious source and target maps. The composition is defined as follows:
$$
(y_0, \gamma \cdot y_0, a)\circ (\gamma \cdot y_0, \delta\gamma \cdot
y_0, b) = (y_0, \delta\gamma y_0,
g(\gamma)\gamma(g(\delta))(g(\gamma\delta))^{-1}ab)
$$
\end{defn}

\begin{thm}
\label{qzg-imag}
The generalised imaginary sort corresponding to the splitting groupoid
of a quantum Zariski geometry is interpretable in the quantum Zariski
geometry.  A quantum Zariski geometry is definable in an algebraically
closed field expanded with the generalised imaginary sort. In
particular, if the splitting groupoid is eliminable, then the quantum
Zariski geometry is interpretable in the pure algebraically closed
field.
\end{thm}

\begin{proof}
  Suppose $V \to X$ is a quantum Zariski geometry. Let $P$ be the set
  of all canonical bases in all fibres $V_x$ up to multiplication, it
  is a definable set and an $S_\bu$-torsor. Indeed, let $\Sigma$ be
  the set of definable maps $h: V \times_X Y \to Y \times \A^n(k)$
  such that $\uu_1, \ldots, \uu_d$ under these maps correspond to
  operators that act by matrices $(t_{ij}^k)$ that satisfy the
  equations $f_l(x,t)=0$ from the clause~\ref{canbase} of the
  Definition~\ref{quantum-zg}. The set $P$ is defined as the set
  $$
  \suchthat{ (v_1, \ldots, v_n ) \in V^n } { \exists h\in \Sigma\,
    \exists y\in Y: h(e_i)=(y,e_i)\ \forall i\ 1  \leq i \leq n }
  $$
  where $\{ e_i \}_{i=1}^n$ is the standard basis in $Y \times
  \A^n(k)$ regarded as a family of vector spaces.  The projection map
  is the restriction of the map $p: V \to X$ to $P$, and the anchor
  map $a: P \to Y$ is defined as
  $$
  a(\bar v) = \suchthat{ y \in Y }{ \exists h \in \Sigma:
    h(e_i)=(y,e_i)\ \forall i\ 1 \leq i \leq n }
  $$
  The action of $S_\bu$ is definable in
  $V^{\eq}$:
  $$
  (y_0, \gamma \cdot y_0, \alpha)\cdot (v_1, \ldots, v_n) = (\alpha \cdot g(\gamma) \cdot
  v_1, \ldots, \alpha \cdot g(\gamma) \cdot v_n)
  $$ 

  Conversely, let $P \to [S_\bu]=X$ be a $S_\bu$-torsor with the
  anchor map $a: P \to S_0$. Let $W=P \times k^n/\G_m$ where $\G_m$
  acts by the formula $a \cdot (p, x) = (a \cdot p, a^{-1} \cdot
  x)$. Define the action of the groupoid $S_\bu$ on $W$ as follows:
  $$
  (y_0, \gamma y_0, a) \cdot (p, x) = ((y_0, \gamma y_0, a) \cdot p,
  g(\gamma)^{-1} x)
  $$
  In particular, $\Gamma$ acts on $W$ if one lets an element $\gamma
  \in \Gamma$ acts on an element $(p,x)$ like the arrow  $(a(p), \gamma
  \cdot a(p), 1)$. Let $V$ be the quotient $W/\Gamma$. The matrices of
  $\uu_k$ define endomorphisms of $W$ which are $\Gamma$-invariant,
  and hence descend to $V$. 
\end{proof}

\begin{defn}[Brauer group]
  Two Azumaya algebras $A,A'$ over $X$ are called \emph{Morita
    equivalent} if there exists a vector bundle $V \to X$ such that $A
  \otimes \End(V) \cong A' \otimes \End(V)$ where $\End(V)$ is the
  $k[X]$-algebra of endomorphisms of the vector bundle $V$. The Morita
  equivalence classes of Azumaya algebras over $X$ form a group under
  the tensor operation called the \emph{Brauer group of $X$}.
\end{defn}

\begin{rem}
  It is not essential in the definition of a quantum Zariski geometry
  that $X$ be affine. One could propose a straightforward
  generalisation of the definition of a quantum Zariski geometry for
  an arbitrary variety using sheaves of
  $\mathcal{O}_X$-algebras. Using such definition it would be possible
  to realise an arbitrary Azumaya algebra in a quantum Zariski
  geometry (in the sense of Propostion~\ref{qzg-surj}).
\end{rem}

\begin{prop}[Theorem~2.5,  \citep{milne}]
  Let $X$ be a variety. Then $\Br(X)$ injects canonically into
  $H^2(X_\et, \G_m)$. 
 \end{prop}

 \begin{prop}[Proposition~2.7, \citep{milne}]
  An equivalence class of an Azumaya algebra of rank $n^2$ specified
  by a 1-cocycle $\alpha \in \check{H}^1(X_\et, \G_m))$ is
  $\delta(\alpha)$ where   $\delta$ is the connecting morphism in the
  following long exact sequence:
  $$
  \ldots \to \check{H}^1(X_\et, \GL_n) \to \check{H}^1(X_\et, \PGL_n)
  \xrightarrow{\delta} \check{H}^2(X_\et, \G_m) \to \ldots
  $$  
 \end{prop}

 \begin{prop}
   \label{torsion}
   Every $n$-torsion element of $H^2(X_\et,\G_m)$ is the image of an
   element of $H^2(X_\et, \mu_n)$.
 \end{prop}

 The fact that the image of the injection is the torsion part of
 $H^2(X_\et, \G_m)$ when $X$ is the spectrum of a field is classical
 (in this case \'etale cohomology coincides with Galois
 cohomology). This fact is established for some other classes of
 varieties, fore example, when $X$ is smooth variety over a field
 (Grothendieck, \citep{grothen-brauer}), when $X$ is affine (due to Gabber,
 unpublished, but see \citep{dejongh}).

 \begin{prop}
   The class of the splitting groupoid in $H^2(k(X), \G_m)$ coincides
   with the restriction of the class of the Azumaya algebra associated
   to a quantum Zariski geometry to the generic point.
 \end{prop}

 \begin{proof}
   The class of the restriction of the splitting groupoid to the
   generic point of $X$ lies in the image of the map
   $H^2(\Gal(k(Y)/k(X)), \G_m) \to H^2(k(X), \G_m)$.  The map $g:
   \Gal(Y/X) \to \PGL_n(Y)$ of Definition~\ref{quantum-zg} is a
   1-cocycle that defines the gluing data for the Azumaya algebra. The
   2-cocycle of the element of an Azumaya algebra associated to the
   quantum Zariski geometry is the coboundary of a the lifting of a
   1-cocycle in $H^1(X, \PGL_n)$ to $\GL_n$ (see \citep{milne},
   Theorem~2.5). From the definition of splitting groupoid and the
   proof of Theorem~\ref{cohomology} it follows that the 2-cocycle
   corresponding to it is $(\sigma,\tau) \mapsto
   g(\sigma)\sigma(g(\tau))(g(\sigma\tau)^{-1})$. One checks that this
   coincides with the definition of the coboundary.
 \end{proof} 
 
 For a smooth variety $X$, the restriction map $\Br(X) \to \Br(k(X))$
 is an injection (Auslander and Goldman, \cite{auslander-goldman}), so
 one can check the eliminability of the splitting groupoid at the
 generic point.

 \begin{rem}
   \label{rem-qzg-groups}
   In view of Propositions~\ref{torsion} and \ref{group-image}, in
   order to eliminate the splitting groupoid of a quantum Zariski
   geometry it suffices to add a generalised imaginary sort for a
   subgroupoid of the splitting groupoid bounded by $\mu_n$.
 \end{rem}
 
\subsection{Quantum Zariski geometries: an example}
\label{sec:quantum-ex}

As an application of the main result of last section I am going to
show that the quantum Zariski geometry corresponding to the quantum
torus algebra (Example~2.1, \citep{quantum-zg}) is not interpretable
in an algebraically closed field.

The input data is as follows. The quantum torus algebra is the algebra 
$$
A = k\langle \uu, \vv, \uu^{-1}, \vv^{-1} \mid \uu\vv=q \vv\uu,
\uu\uu^{-1}=\uu^{-1}\uu=1, \vv\vv^{-1}=\vv^{-1}\vv=1 \rangle 
$$
where $q^n=1$ for some integer $n > 1$. The cover $Y \to X$ is defined
by the equations $\mu^n = x, \nu^n = y$, where $x,y$ are the coorinate
function on $X = \G_m \times \G_m$. The matrices of $\uu, \vv$ are
$$
\uu = \left(
  \begin{array}{cccc}
    \mu & 0 & \ldots & 0\\
    0 & q \mu  & \ldots & 0\\
    \vdots &  & & \vdots \\
    0  &  & \ldots & q^{n-1} \mu \\
  \end{array}\right)
\qquad\qquad \vv = \left(
  \begin{array}{ccccc}
    0 & 0 & \ldots & 0 & \nu\\
    \nu & 0 & \ldots & 0 & 0\\
    0 & \nu & \ldots & 0 & 0\\
    \vdots & &  & & \vdots \\
    0 & &  \ldots & \nu & 0\\
  \end{array}\right)
% \nu e_{(i+1) \mod n},\qquad 0 \leq i < n
$$
and the matrices of $\uu^{-1}, \vv^{-1}$ are inverses of these.  The
Galois group of the cover is the product of two cyclic groups of order
$n$, and the cocycle $g$ is specified on generators $\alpha, \beta$ of
cyclic factors of $\Gamma$ by matrices
$$
g(\alpha)=\left(
  \begin{array}{ccccc}
   0 & 0 & \ldots & 0 & 1\\ 
   1 & 0 & \ldots & 0 & 0\\ 
   0 & 1 & \ldots & 0 & 0\\ 
   \vdots &   &   &   & \vdots\\ 
   0 & 0 & \ldots & 1 & 0\\ 
  \end{array}
\right) \qquad
g(\beta)=\left(
  \begin{array}{ccccc}
   1 & 0 & 0 & \ldots & 0\\ 
   0 & \zeta & 0 & \ldots & 0\\ 
   \vdots &   &   &   & \vdots\\ 
   0 & 0 & 0 & \ldots & \zeta^{n-1} \\ 
  \end{array}
  \right)
$$

\begin{defn}
\label{defn:cyclic-csa}
Let $k$ be a field that contains $n$-th root of unity and such that
$n$ is invertible in $k^\times$. A \emph{cyclic algebra of rank $n$}
corresponding to elements $a, b \in k^\times$ is the algebra
$$
k\langle x, y \mid x^n = a, y^n = b, xy=\zeta yx \rangle
$$
where $\zeta$ is a primitive $n$-the root of unity.
\end{defn}

\begin{prop}[\citet{gille-csa}, Corollary~4.7.4]
\label{prop:cyclic-csa}
The class of a cyclic algebra that corresponds to an element $b \in
k^\times$ is trivial in $Br(k)$ if and only if $b \in
N_{K/k}(K^\times)$ where $K = k(\sqrt[n]{a})$ and
$N_{K/k}(x)=\prod_{\sigma \in \Gal(K/k)} \sigma(x)$.
\end{prop}

The restriction of the Azumaya algebra of the quantum torus Zariski
structure to the generic point is clearly a cyclic algebra
corresponding to coordinate functions $x,y$ of $X$. In order to check
if this algebra is split one has to verify if $y \in
N_{K(\mu)/K}(K(\mu)$ where $K=k(X)$ which is clearly false as for $f
\in K(\mu)$ $N_{K(\mu)/K}(f(\mu))$ only has terms with coefficients
that contain $y$ to the power that divides $|\Gal(Y/X)|$.

Let us now illustrate the point made in
Remark~\ref{rem-qzg-groups}. The cohomology class $h \in H^2(X, \G_m)$
that corresponds to the quantum torus structure is the image of a
cocycle in $H^2(X, \G_m \times \G_m)$ which corresponds to the central
extension of the Galois group of the cover $Y \to X, (\mu, \nu)
\mapsto (\mu^n, \nu^n)$ by $\mu_n$:
$$
1 \to \mu_n \to G \to \Z/n\Z \times \Z/n\Z \to 1
$$
where $G=\langle u,v \mid u^n=v^n=1, [u,[u,v]]=[v,[u,v]]=1,
[u,v]^n=1\rangle$. Let us interpret the quantum Zariski geometry in
the group extension Zariski geometry $M$ associated to this short
exact sequence.

Following the construction of Theorem~\ref{qzg-imag}, let $W$ be the
quotient $M \times k^n/\mu_n$ and let $V$ be the quotient of $W$ by an
action of $G$ (which is the same as the action of an action groupoid)
defined as follows. The action of $G$ on $M$ is already defined, so we
have to define the action on $k^n$ depending the point of $y$ . Let
the action of an element $\gamma \in G$ be given by
$g(\gamma/\mu_n)$. Define the action of $\uu, \vv$ to be the action of
the generators $u, v$ of the group $G$ on $V$. This defines an
interpretation of the quantum torus geometry in $M$. 
\subsection{Eliminability of generalised imaginaries in \textrm{CCM}
  and \textrm{RCF}} 
\label{sec:elim-ccm-omin}

We have seen that generalised imaginaries serve as an obstruction to
intrerpetability of quantum Zariski geometries in the theory of
algebraically closed fields. One easily observes that groupoids
definable in the projective sort of the compact complex manifolds
structure are not any more eliminable in this structure than they are
in algebraically closed fields.

\begin{prop}
  Let $X_\bu$ be a groupoid definable in
  the field $\C$. Consider the same groupoid as a definable groupoid
  in the projective sort of the compact complex spaces structure, call
  it $Y_\bu$. Then $X_\bu$ is eliminable if and only if $Y_\bu$ is
  eliminable.
\end{prop}

\begin{proof}
  By Proposition~\ref{finite-quot} any groupoid is Morita equivalent
  to a groupoid $X_\bu$ such that $X_0 \to [X_\bu]$ has finite
  fibres. Their eliminability depends on the Galois groups of generic
  points of definable subsets of $X_0$ (Theorem~\ref{cohomology}).

  The statement of the proposition follows from purity of the field
  $\C$ intepreted in the projective sort in the CCM. As a
  consequence all Galois groups of sets of parameters in $\C \subset
  \mathbb{P}^1$ coincide with that of algebraically closed field, and
  eliminability does not change.
\end{proof}

On the other hand, $ACF$-definable groupoids are always eliminable in
$RCF$. 

\begin{prop}
  Let $X_\bu$ be a groupoid bounded by an Abelian group definable in
  an algebraically closed field $k$ over a set of parameters $K$ and
  let $X_\bu$ be definable in $R$ in such a way that the structure
  maps are definable in the field $k$ identified with $R +
  \sqrt{-1}R$.  Then $X_\bu$ is eliminable in $R$ (over the same set
  of parameters $K$).
\end{prop}

\begin{proof}
  Let $p: X_0 \to [X_\bu]$ be the projection on the definable set of
  connected components. Since $p$ has a section $j$ definable in $R$
  there exists a definable $X_\bu$-torsor over $[X_\bu]$: $Q = \cup_{x
    \in X_0, y \in \mathrm{Im}(j)} \Mor(x,y)$ with the action of $X_1$
    given by composition of arrows.
\end{proof}

In particular, quantum Zariski geometries are interpetable in
$\R$, and they are not interetabile in the structure of compact
complex spaces any more than they are interpretable in an
algebraically closed field.

\appendix

\renewcommand{\theequation}{\Alph{section}.\arabic{equation}}
\numberwithin{thm}{section}

\section{Group cohomology}
\label{groupcoho}

In this appendix the necessary facts about group cohomology are
recalled, for detailed exposition see \citep{serre-galois, neukirch}.

\begin{defn}[Group cohomology]
  Let $G$ be a profinite group. A \emph{$G$-group} $A$ is a disceret
  group that is endowed with a continuous action of $G$, i.e. a
  continuous homomorphism $G \to \Aut(A)$. If $A$ is Abelian then $A$
  is called a \emph{$G$-module}.

  Let $A$ be a $G$-module. The group cohomology functors $H^i(G,A)$
  are the derived functors of the functor $(-)^G$ that takes a
  $G$-module $A$ to its submodule of $G$-invariant elements: $A^G = \{
  a \in A \mid \forall g \in G\ ga=a \}$.
\end{defn}

Any $G$-module has a standard acyclic resolution that gives rise
to the \emph{homogeneous} cochain complex that computes the
cohomology. This complex is quasi-isomorphic to the
\emph{inhomogeneous} cochain complex. Computing its cohomology of the
latter complex the following elementary definition of the functors
$H^i(G, \text{--})$.

\begin{defn}
  Let $A$ be a $G$-module. The $n$-th term of the inhomogeneous cochain
  complex $C^n(G,A)$ is defined to be the set of all continuous maps
  $G^n \to A$ with $C^n(G,A)$ supposed formally to be 0 for $n<0$ and
  $C^0=A$. The differential $d_n: C^n(G,A) \to C^{n+1}(G,A)$ is
  defined as follows
  \begin{eqnarray*}
    (dh)(\sigma_0, \ldots, \sigma_n) & = & \sigma_1 \cdot
    h(\sigma_1, \ldots, \sigma_n) + \\
    & & + \sum_{i=1}^nh(\sigma_0, \ldots,
    \sigma_{i-1}\sigma_i, \ldots, \sigma_n) + (-1)^nh(\sigma_0,
    \ldots, \sigma_{n-1})
  \end{eqnarray*}
  The $n$-cohomology group of the complex, $\Ker d_n/\Im d_{n-1}$ is
  called the $n$-th cohomology group of $G$ with coefficients in the
  module $A$, and is denoted $H^n(G,A)$.
  In particular, the 1-cocycles are maps $h: G \to A$ such that 
  $$
  h(\sigma\tau) = h(\sigma) + \sigma\cdot h(\tau)
  $$
  modulo the equivalence relation: $h \sim h'$ if and only if there
  exists $g \in A$ such that $h = \sigma(g) + h' - g$. This definition also makes sense
  when $A$ is non-Abelian, though $H^1$ has no group structure and is
  just a set with a distinguished element of cocycles cohomolougous to
  the zero cochain.

  The second cohomology group is the set of maps $G^2 \to A$ such that
  $$
  h(\alpha\sigma,\tau) = h(\alpha,\sigma\tau) - h(\alpha,\sigma) +
  \alpha\cdot h(\sigma,\tau)
  $$
  and two 2-cocycles $h,h'$ are cohomologous if there exists a
  map $g: G \to A$ such that 
  $$
  h(\sigma,\tau) = h'(\sigma,\tau) + g(\sigma) - g(\sigma\tau) +
  \sigma\cdot g(\tau)
  $$
\end{defn}

The cohomology of profinite group is related to the cohomology of
element of the inverse systems in the expected way.

\begin{prop}[\citep{serre-galois}, Chapter~I,Corollary~2.2]
\label{coho-limit}
  Let $G$ be a profinite group and let $A$ be a $G$-module. Then
  $$
  H^i(G,A) = \varprojlim_{U \subset G} H^i(G/U, A^U)
  $$
  where the limit is taken over all closed subgroups $U \subset G$.
\end{prop}

The low degree cohomology groups have natural geometric and
group-theoretic interpretations.

Let $A$ be an Abelian algebraic group. Then the set of $A$-torsors
over $K$ which have a point in an extension $L$ forms an Abelian group
called \emph{Weil-Ch\^atelet} group, which is isomorphic to
$H^1(\Gal(L/K), A)$. The group operation is defined as follows. Let
$P,Q$ be two $A$-torsors. Then the result of the group operation is
the quotient of $P \times Q$ by the action of $A$: $a\cdot (p,q) = (a
\cdot p, a^{-1}\cdot q)$ (note the similarity with the Baer sum of
group extensions or groupoids). If $A$ is 0-dimensional then the
existence of the quotient is straightforward, while in general it
takes some work to construct it (which involves Weil's theorem on
birational group laws, see \citep{weil} for details).

\begin{thm}[Kummer theory]
  Let $K$ be a field that contains $n$-th roots of unity with $n
  \not| \chara K$. Then $H^1(K,
  \mu_n) = K^\times/(K^\times)^n$.
\end{thm}

Let $A$ be an Abelian group and let $1 \to A \to G \xrightarrow{p} H
\to 1$ be a group extension. Regard $A$ as an $H$-module: an element
$b\in H$ acts on $a$ by conjugation, $b\cdot a = j(b)\cdot a \cdot
j(b)^{-1}$, for some section $j$ of $p$, the action is independent of
the choice of section since $A$ is normal. Now fix a section and
associate to it a cocycle:
$$
f(\sigma,\tau) = j(\sigma) \cdot j(\tau) \cdot j(\sigma\tau)^{-1}
$$
Conversely, given a cocycle $f$ one defines a group law on $A \times
H$
$$
(a,\sigma) \cdot (b, \tau) = (a + \sigma\cdot b + f(\sigma,\tau),
\sigma\tau)
$$
\begin{thm}[\citet{neukirch}, Theorem~1.2.5]
  \label{group-ext}
  There is a bijective correspondence between elements of $H^2(H,A)$
  and extensions of $H$ by $A$ such that $A$ has the prescribed
  $H$-module structure.
\end{thm}
Note that if $A$ is a definable group and $H$ is finite, the group
extension is also a definable group.

\begin{prop}
  Let $A$ be an Abelian group and let $1 \to A \to G \xrightarrow{p}
  H$ be a split group extension. Then the set of all sections of $p: G
  \to H$ modulo conjugation by elements of $H$ is a torsor under
  $H^1(H,A)$: given a section $j$ put $(h \cdot j)(a) =
  j(a)(h(ab))^{-1}$.
\end{prop}
Indeed, the cocycle condition tells us how $h(a)$ commutes with $j(a)$
and
$$j(ab)h(ab)=j(a)j(b)(j(b))^{-1}h(a)j(b)h(b)=j(a)h(a)j(b)h(b)$$ so
$(h \cdot j)$ is a homomorphism; similarly, one checks that acting by
coboundaries conjugates a section by an element of $H$.

\bibliography{imag}

\noindent {\sc Einstein Institute of Mathematics \\
Edmond J. Safra Campus\\
The Hebrew University of Jerusalem\\
Givat Ram. Jerusalem, 91904, Israel\\}
{\tt sustretov@ma.huji.ac.il\\
}

\end{document}